\documentclass[titlepage,twoside,12pt]{article}
\usepackage{amsmath}
\usepackage{amssymb}
\usepackage{latexsym}
\usepackage{graphicx}
\usepackage{amsfonts}
\begin{document}

{\LARGE \bf Brief Lecture Notes on Self-Referential \\ \\ Mathematics, and Beyond} \\ \\

{\bf Elem\'{e}r E Rosinger} \\ \\
Department of Mathematics \\
and Applied Mathematics \\
University of Pretoria \\
Pretoria \\
0002 South Africa \\
eerosinger@hotmail.com \\ \\

{\bf Abstract} \\

Recently delivered lectures on Self-Referential Mathematics, [2], at the Department of Mathematics and Applied Mathematics, University of Pretoria, are briefly presented. Comments follow on the subject, as well as on Inconsistent Mathematics. \\ \\

{\bf 0. Prologue} \\

The basic idea in the Self-Referential Mathematics, [2], is to replace the Foundation Axiom, (FA), in Set Theory with a
suitable Anti-Foundation Axiom, (AFA), in such a way as to : \\

- keep all the sets in the usual Set Theory \\

while at the same time, to : \\

- allow a large class of new sets, sets given this time by {\it self-referential} definitions. \\

In other words, one obtains a significant extension of usual Set Theory, an extension which is proved to be consistent,
provided usual Set Theory is consistent. \\

As it happens not seldom in science, the terminology used may turnout to be rather inappropriate, if not in fact misleading. The same
happens in [2], where the term "vicious circle" is used instead of "self-referential". \\

A likely reason for that particular terminology in [2], which has a clear negative connotation, comes from the fact that the 1903
Russell Paradox in Set Theory is based on self-reference, being but a reformulation in mathematical, in particular, set theoretic terms of the ancient
Greek paradox of the liar. \\

On the other hand, when considered in the larger and longer perspective of human tradition and civilization, self-reference, together with infinity and
change, have since the earliest known, in fact, prehistoric times been some of the fundamental ideas preoccupying human thought, and as such, they have
not had any sort of inevitable negative connotation, see section 5. It follows, therefore, that the term "vicious" can be seen as an overstatement
resulted from a partial view of what self-reference does in fact encompass and mean in its more full generality. \\

We can in essence clarify as follows the aims and the means of the Self-Referential Mathematics in [2]. Let us consider the following three groups of
axioms of Set Theory, see section 6 for all the usual axioms, in particular, those used in [2] : \\

ZFC$^-$ = Zermello - Fraenkel + Choice \\

ZFC = ZFC$^-$ + FA \\

ZFA = ZFC$^-$ + AFA \\

where the AFA axiom with be specified in section 2. \\

At first, it may appear that the Set Theories corresponding to ZFC and ZFA may be rather different, since their common part corresponds only to ZFC$^-$,
while the respective additional axioms FA and AFA seem in fact to be inconsistent with one another. \\
However, as it turns out this is not the case. And what happens instead is that :

\begin{itemize}

\item The Set Theory based on the ZFA axioms contains all the sets in the Set Theory based on the ZFC, and in addition, contains a large class of other
sets obtained by self-referential definitions. \\

\item The axioms ZFA are consistent, provided that the axioms ZFC are consistent.

\end{itemize}

As for the traditional and still exclusively predominant idea of the absolute necessity of consistency, one should
consider the recent emergence of Inconsistent Mathematics, see [11,12]. And in fact, as far as everyday practice is concerned,
we have for more than half a century by now been basing much of our lives on a specific form of Inconsistent Mathematics.
Indeed, our ever more pervasive and critically important electronic digital computers are - even when only operating on
non-negative integers - functioning according to the Peano Axioms, plus the Machine Infinity Axiom, namely \\

$\exists~~ M >> 1 ~:~ M + M = 1$ \\

where $M$ is called "machine infinity", and typically is larger than $10^{1000}$. And obviously, the Peano Axioms are
trivially inconsistent with the Machine Infinity Axiom. \\ \\

{\bf 1. Sets, Ur-Elements and Classes} \\

We denote by SET the {\it class} of all sets, including the unique void set $\phi$. As is well known, a lot of mathematics can be built up starting alone
with the void set $\phi$. Indeed, as a first step, and following von Neumann, one can define the non-negative integers by \\

$ 0 = \phi,~ 1 = \{ \phi \},~ 2 = \{ \{ \phi \} \}, \ldots $ \\

and then, step by step build all the integers, the rational and real numbers, and so on. Further, one can define Cartesian products, binary relations,
functions, etc., and obtain a considerable part of mathematics in this manner. \\

In the sequel, it will be convenient to allow, in addition to the void set $\phi$, other such starting entities in the construction of mathematics. The
class of such entities is denoted by ${\cal U}$, and any respective element $u \in {\cal U}$ is called an {\it ur-element}, assumed to have only one
property, similar to that of the void set $\phi$, namely that the relation \\

$~~~~~~~~~~~~ a \in u$ \\

does {\it not} hold for any entity $a$ in the theory. \\

In this way, there will be three types of entities in the theory, namely \\

1)~~~ $SET$, which is the class of all sets, \\

2)~~~ ${\cal U}$ which is the class of all ur-elements, and \\

3)~~~ $CLASS$ which denotes all the classes. \\

Here it is understood that any set $a \in SET$ is a "small" class, while $SET$ itself is one of the "proper" classes, since it is not itself a set. In
other words, $SET \in CLASS$, $SET \notin SET$, $CLASS \notin SET$. \\
As for ur-elements, it is assumed that ${\cal U} \in CLASS$ and it is another instance of "proper" class, thus in particular ${\cal U} \notin SET$. \\

Briefly, we have therefore \\

1)~~~ $SET$ denotes all the sets, and it is a proper class \\

2)~~~ $CLASS$ denotes all the classes \\

3)~~~ ${\cal U}$ denotes all the ur-elements, and it is a proper class \\

4)~~~ a set is a "small" class \\

5)~~~ a class which is not a set is "large", thus it is a {\it proper} class \\

6)~~~ an ur-element does not have any elements, either sets, ur-elements, \\
      \hspace*{0.65cm} or classes \\

7)~~~ non-set = proper class $\bigvee$ ur-element \\

8)~~~ non-set $\bigwedge$ non-class = ur-element \\

9)~~~ every predicate determines a class \\

10)~~~ a subclass of a set is a set \\

11)~~~ sets are closed under a number of operations, among them, \\
       \hspace*{1cm} pairing, union, power set, see below \\

12)~~~ $a \in b \in CLASS ~~\Longrightarrow~~ a \in SET$ \\

13)~~~ the class SET of all sets is "large", thus it is a proper class \\

\hspace*{2cm} indeed, according to Russell's Paradox, let SET be a set, \\
\hspace*{2cm} then $R = \{ a \in SET ~|~ a \notin a \}$ is a set, thus $R \in SET$, and
\hspace*{2cm} therefore $R \in R ~~\Longleftrightarrow~~ R \notin R$, which is absurd \\

As for the binary relation $\in$, we have \\

\hspace*{0.5cm} - set $\in$ set \\
\hspace*{0.5cm} - set $\in$ proper class \\
\hspace*{0.5cm} - set $\notin$ ur-element \\
\hspace*{0.5cm} - proper class $\notin$ set \\
\hspace*{0.5cm} - proper class $\notin$ proper class \\
\hspace*{0.5cm} - proper class $\notin$ ur-element \\
\hspace*{0.5cm} - ur-element $\in$ set \\
\hspace*{0.5cm} - ur-element $\in$ proper class \\
\hspace*{0.5cm} - ur-element $\notin$ ur-element \\

thus denoting $\in$ by $\rightarrow$, while $\notin$ by $\nrightarrow$ , we have \\ \\

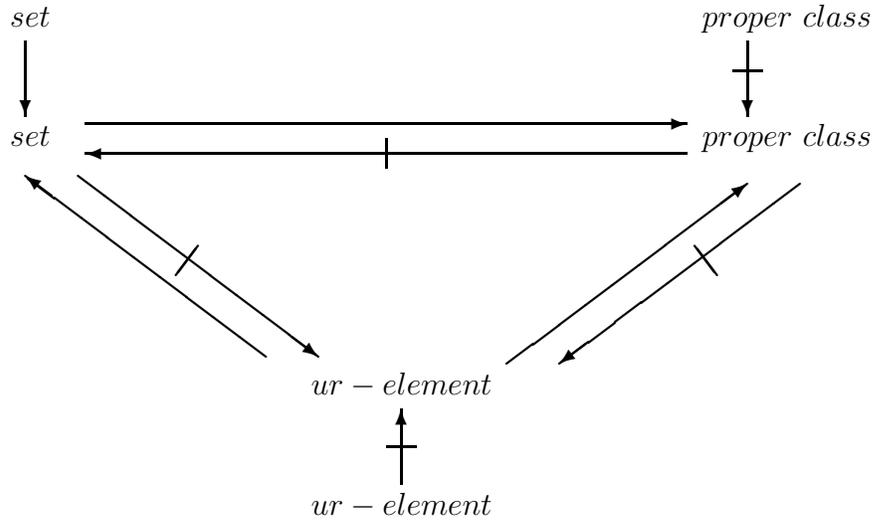
\begin{math}
\setlength{\unitlength}{0.2cm}
\thicklines
\begin{picture}(60,29)

\put(5,28){$set$}
\put(6,27){\vector(0,-1){5}}
\put(5,20){$set$}
\put(10,21.5){\vector(1,0){40}}
\put(50,19.5){\vector(-1,0){40}}
\put(30,18.5){\line(0,1){2}}
\put(51,28){$proper~ class$}
\put(54,27){\vector(0,-1){5}}
\put(53,25){\line(1,0){2}}
\put(51,20){$proper~ class$}
\put(9.5,18){\vector(4,-3){16}}
\put(16,11.4){\line(3,4){1.5}}
\put(22,6){\vector(-4,3){16}}
\put(25,3.5){$ur-element$}
\put(25,-4.5){$ur-element$}
\put(31,-2.5){\vector(0,1){5}}
\put(30,0){\line(1,0){2}}
\put(38,5.5){\vector(4,3){16}}
\put(57.5,17.5){\vector(-4,-3){16}}
\put(52,11.4){\line(-3,4){1.5}}

\end{picture}
\end{math} \\ \\ \\ \\

\begin{math}
\setlength{\unitlength}{0.2cm}
\thicklines
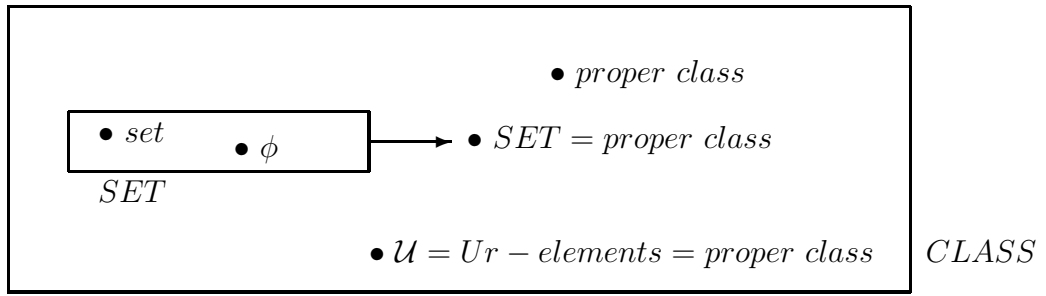
\begin{picture}(60,21)

\put(1,20){\line(1,0){60}}
\put(1,20){\line(0,-1){19}}
\put(1,1){\line(1,0){60}}
\put(61,20){\line(0,-1){19}}

\put(5,13){\line(1,0){20}}
\put(5,13){\line(0,-1){4}}
\put(5,9){\line(1,0){20}}
\put(25,13){\line(0,-1){4}}

\put(7,11){$\bullet~ set$}
\put(16,10){$\bullet~ \phi$}
\put(37,15){$\bullet~ proper~ class$}
\put(7,7){$SET$}
\put(31.5,10.6){$\bullet~ SET = proper~ class$}
\put(25,11){\vector(3,0){5.5}}
\put(25,3){$\bullet~ {\cal U} = Ur-elements = proper~ class$}
\put(62,3){$CLASS$}

\end{picture}
\end{math} \\

{\bf Note on ur-elemets} \\

Clearly \\

$ a \notin \phi$ \\

holds for all entities $a$ in the theory, however, it is nevertheless considered that \\

$ \phi \in SET$ and $\phi \notin {\cal U}$ \\

Also, it is possible that \\

$ X \in SET$ and $X \cap {\cal U} \neq \phi$, or even $X \subseteq {\cal U}$ \\

For instance, if $x \in {\cal U}$, then \\

$ X = \{ x \} \in SET,~~~ X \subseteq {\cal U}$. \\

{\bf Note on $SET$} \\

Here we should clarify that $SET$ denotes, in fact, all the sets which exist in the Set Theory based on the ZFA axioms. Therefore, let us denote by
$SET_0$ all the sets in the Set Theory based on the ZFC axioms. Then as mentioned in section 0, and seen later, we have $SET_0 \subsetneqq SET$, thus the
above diagram can be made more precise as follows \\ \\

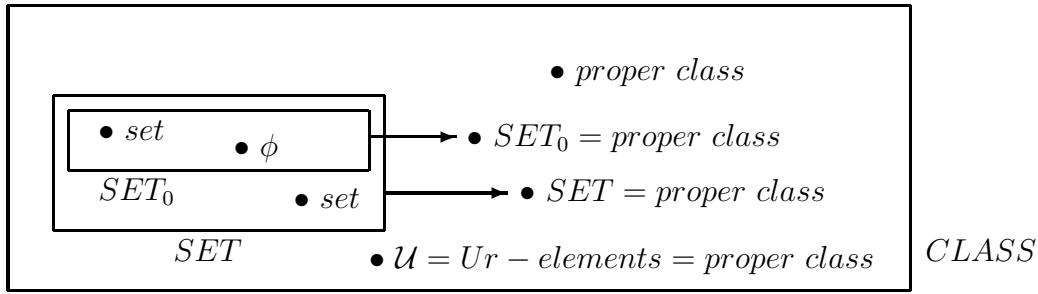
\begin{math}
\setlength{\unitlength}{0.2cm}
\thicklines
\begin{picture}(60,21)

\put(1,20){\line(1,0){60}}
\put(1,20){\line(0,-1){19}}
\put(1,1){\line(1,0){60}}
\put(61,20){\line(0,-1){19}}

\put(4,14){\line(1,0){22}}
\put(4,14){\line(0,-1){9}}
\put(4,5){\line(1,0){22}}
\put(26,14){\line(0,-1){9}}

\put(5,13){\line(1,0){20}}
\put(5,13){\line(0,-1){4}}
\put(5,9){\line(1,0){20}}
\put(25,13){\line(0,-1){4}}

\put(7,11){$\bullet~ set$}
\put(16,10){$\bullet~ \phi$}
\put(37,15){$\bullet~ proper~ class$}

\put(12,3){$SET$}
\put(20,6.5){$\bullet~ set$}
\put(35,7){$\bullet~ SET = proper~ class$}
\put(26,7.5){\vector(3,0){8.2}}

\put(7,7){$SET_0$}
\put(31.5,10.6){$\bullet~ SET_0 = proper~ class$}
\put(25,11.25){\vector(3,0){5.8}}
\put(25,2.5){$\bullet~ {\cal U} = Ur-elements = proper~ class$}
\put(62,3){$CLASS$}

\end{picture}
\end{math} \\

{\bf Examples of Sets} \\

There are only {\it two} kind of sets $a \in SET$, namely \\

(1.1)~~~ $a = \phi$, which is equivalent with $\neg \,(~ \exists~~ b \in SET \bigcup {\cal U} ~:~ b \in a ~)$ \\

(1.2)~~~ $\exists~~ b \in SET \, \bigcup \, {\cal U} ~:~ b \in a$ \\

{\bf Operations with sets} \\

An {\it ordered pair} is the set ~$< a, b > ~=~ \{ \{ a \}, \{ a, b \} \}$, with $a, b \in$ SET \\

thus ~$< a, b > ~=~ < c, d > ~~\Longleftrightarrow$~~ a = c,~ b = d \\

A {\it relation} $R \in SET$ has all its elements given by pairs ~$< a, b >$, where $a \in A,~ b\in B$, for two suitably
given sets $A,~ B \in SET$. Often for convenience one denotes $a R b$ for ~$< a, b > \, \in R$. \\

If $A \in SET$ is such that ~$< a, b >\, \in R ~\Longrightarrow~ a, b \in A$, then $R$ is a relation on $A$. \\

A {\it relational structure} is ~$< A, R >$,~ with $R$ relation on $A$. \\

A {\it function} is a relation $R$ such that ~$< a, b >,\, < a, c >\, \in R ~~\Longrightarrow~~ b = c$ \\

If ~$f$~ is a function, then \\

(1.3)~~~ $ dom ( f ) = \{ a ~|~ \exists~~ b ~:~ f ( a ) = b \} $ \\

(1.4)~~~ $ rng ( f ) = \{ b ~|~ \exists~~ a ~:~ f ( a ) = b \} $ \\

thus \\

(1.5)~~~ $f \in c \to d ~~\Longleftrightarrow~~ c = dom ( f ),~ rng ( f ) \subseteq d $ \\

The {\it power set} of $a \in SET$ is \\

(1.6)~~~ ${\cal P} ( a ) = \{ b ~|~ b \subseteq a \}$ \\

Example : if $a = \{ \phi, p \}$, with $p \in SET \bigcup {\cal U}$, then ${\cal P} ( a ) =
\{ \phi, \{ \phi \}, \{ p \}, a \}$. \\

Consider the {\it predicate} P ( x ) given by \\

(1.7)~~~ $x$ is an ordered pair $< a, b >$\, and $b = {\cal P} ( a )$ \\

then this defines the {\it power set} function ${\cal P} : SET \longrightarrow SET$, and \\

(1.8)~~~ $\bigcup dom ( {\cal P} ) = SET$ \\

thus it is "large", and therefore, a proper class. \\

The {\it natural numbers} are \\

(1.9)~~~ $0 = \phi,~ 1 = \{ 0 \} = \{ \phi \},~ 2 = \{ 0, 1 \} = \{ \phi, \{ \phi \} \}, \ldots$ \\

{\it Disjoint union} is $A + B = ( \{ 0 \} \times A ) \cup ( \{ 1 \} \times B )$ \\

For $a \in SET$, we define \\

(1.10)~~~ $\bigcup a = \{~ x ~|~ \exists~ y \in a ~:~ x \in y ~\} = \{ x \in y \in a \}$ \\ \\

{\bf 2. Systems of Equations Which Define Sets} \\

In usual, that is, ZFC Set Theory, one way to define a set $X$ is by an equation \\

$X = \{ x ~|~ P ( x ) \}$ \\

where $P$ is a suitable predicate. Within ZFC, an essential restriction on $P$ is that it {\it cannot} in any way refer to
the set $X$ which it is supposed to define. This condition is meant to avoid a "vicious circle", or in more proper terms,
self-referentiality, an avoidance which has until recently been universally accepted, and in fact required, since
Russell's paradox. \\

In particular, one cannot define any set $a \in SET_0$, even by such a simple equation, like \\

(2.1)~~~ $ a = \{ a \} $ \\

since obviously, it is a self-referential equation. On the other hand, as seen in 1), 4), 5) in Examples 2.1. below, this
equation can easily be solved in $SET$, that is, based on the Anti-Foundation Axiom, (AFA). \\

Here we can note that one cannot define any set $a \in SET_0$, or for that matter, $a \in SET$, even by the yet more
simple equation \\

(2.2)~~~ $ a = a $ \\

since this equation will obviously {\it not} give a unique set in $SET_0$, or in $SET$. \\

Also, as seen in 3) in Proposition 2.3. below, one cannot define a set $a \in SET$ by the equation \\

(2.3)~~~ $ a = {\cal P} ( a )$ \\

Let us consider now the equation \\

$~~~~~~ x = \{ a, x  \}$ \\

where $a \in SET \bigcup {\cal U}$ is given. Then \\

$~~~~~~  x = \{ a, x  \} = \{ a, \{ a, x  \} \} = \{ a, \{ a, \{ a, x \} \} \} = \ldots$ \\

thus an intuitive solution would be \\

$~~~~~~ x = \{ a, \{ a, \{ a, \ldots \} \} \}$ \\

which however is not possible within ZFC, since it would obviously lead to the {\it infinite} descending sequence \\

$~~~~~~ \ldots x \in x \in x \in x \in x$ \\

thus contradict the Foundation Axiom, (FA), see below. \\

Let us now return to the general situation and enquire what should the solution given by sets be of a corresponding system of equations. Let us as an
example consider for that purpose the following system of equations, where $p, q \in SET \bigcup {\cal U}$ are given, and where we want to find sets $x,
y, z \in SET$, such that \\

{}~~~ $ x = \{ x, y \} $ \\

{}~~~ $ y = \{ p, q, y, z \} $ \\

{}~~~ $ z = \{ p, x, y \} $ \\

Let $e : X = \{ x, y, z \} \longrightarrow $ the right-hand sides of the above equations \\

thus \\

{}~~~ $ e_x = \{ x, y \},~ e_y = \{ p, q, y, z \},~ e_z = \{ p, x, y \} $ \\

What is then a solution $s$ to these equations supposed to be ? \\

One way is given by $s : X \longrightarrow SET$, namely $ X \ni v \longmapsto s_v \in SET$, with \\

{}~~~ $ s_x = \{ s_x, s_y \},~ s_y = \{ p, q, s_y, s_z \},~ s_z = \{ p, s_x, s_y \} $ \\

or equivalently \\

{}~~~ $ \begin{array}{l}
             \forall~ v \in X ~: \\ \\
             ~~~~ s_v = \{ s_w ~|~ w \in e_v \cap X \} \, \bigcup \, \{ w ~|~ w \in e_v \cap A \} = \\ \\
             ~~~~~~~ = s [ e_v \cap X ] \, \bigcup \, ( e_v \cap A )
         \end{array} $ \\ \\

where $A = \{ p, q \}$ \\

Returning now to the equation, see (2.2) \\

{}~~~ $x = x$ \\

one way to avoid the inconvenience of non-unique solutions in $SET$, is to take, see also (2.10) below \\

{}~~~ $X = \{ x \} \subseteq {\cal U}$ \\

and then the solution $s$, if it exists, is a {\it function} \\

{}~~~ $s : X \ni x \longrightarrow s_x \in SET$ \\

This liberty to {\it distinguish} between {\it indeterminates}, and on the other hand, the {\it solution} is in fact familiar from usual algebra. Indeed,
if for instance we have the system of equations in real numbers \\

{~}~~~ $2 x + 3 y = 7 \\
{~}~~~~ 5 x - 4 y = 1$ \\

then the set of indeterminates is $X = \{ x, y \}$, while a solution $s$, which in this case exists, is given by a function $s : X \ni v \longrightarrow
s_v \in \mathbb{R}$. This distinction is even more obvious when, as with a system of equations like \\

{~}~~~ $3 x + 2 y = 5 \\
{~}~~~~ 5 x - 3 y = 2$ \\

the indeterminates have the same value $x = y = 1$, thus the solution cannot be identified with the single number $1$, but only with the function $s :
X \ni v \longrightarrow s_v \in \mathbb{R}$, for which $s_x = 1,~ s_y = 1$. \\

Finally we can note that, as seen in Definition 2.2. below, the requirement $X \subseteq {\cal U}$ can on occasion be done away
with. \\

Now, the approach in [2] to Self-Referential Mathematics will be able to accept for sets in $SET$ such definitions which are
given by certain systems of equations that can be {\it self-referential}. There are in this regard {\it three}
kind of systems of equations considered so far. The first kind of systems is given by \\

{\bf Definition 2.1.} \\

A structure ${\cal E} ~=~ < X, A, e >$ is called a {\it flat system of equations}, if and only if \\

(2.4)~~~ $ X,~ A \in SET$,~~~ $X \subseteq {\cal U},~~~ X \bigcap A = \phi $ \\

with $X$ the set of $indeterminates$ and $A$ the set of $atoms$, while \\

(2.5)~~~ $ e : X \longrightarrow {\cal P} ( X \bigcup A ) $ \\

defines the $right~ hand~ terms$ of the equations of the system, see (2.8) below, with \\

(2.6)~~~ $ X \ni v \longmapsto b_v = e_v \cap X $ \\

being the set of indeterminates on which $v$ immediately depends, and similarly, with \\

(2.7)~~~ $ X \ni v \longmapsto c_v = e_v \cap A $ \\

being the set of atoms on which $v$ immediately depends \\

In this way, the flat system of equations is given by \\

(2.8)~~~ $ x = e_x,~~~ x \in X $ \\

which of course can in particular be one single equation, when $X = \{ x \}$ is a set with one single element. \\

A $solution$ to ${\cal E}$ is a function \\

(2.9)~~~ $ X \ni x \longmapsto s_x = \{ s_y ~|~ y \in b_x \} \bigcup c_x \in SET $ \\

and we denote \\

$solution-set ({\cal E}) = \bigcup \, \{ s_x ~|~ x \in X \} =$ \\

$~~~~~~ =  \{ s_y ~|~ y \in e_x \cap X,~ x \in X \} \, \bigcup \, \bigcup_{x \in X} ( e_x \cap A ) = s [ X ] \in SET $ \\

as well as \\

$V [ A ] = \bigcup \left \{~ solution-set ({\cal E}) ~ \begin{array}{|l}
                                                         ~{\cal E} ~\mbox{flat system of equations} \\
                                                         ~\mbox{with atoms}~ A
                                                      \end{array} ~\right\} = $ \\ \\

$~~~~~~ = \left \{ c ~~ \begin{array}{|l}
                                      ~\exists~ {\cal E} =\, < X, A, e > ~\mbox{flat system of equations} : \\
                                      ~~~~ c \in solution-set ({\cal E})
                                 \end{array} \right\} = $ \\ \\

$~~~~~~ = \left \{ s_y ~~ \begin{array}{|l}
                                      ~\exists~ {\cal E} =\, < X, A, e > ~\mbox{flat system of equations} : \\
                                      ~~~~ y \in e_x \cap X,~~\ x \in X
                                 \end{array} \right\} \, \bigcup $ \\ \\

$~~~~~~ \bigcup \, \left \{ c ~~ \begin{array}{|l}
                                      ~\exists~ {\cal E} =\, < X, A, e > ~\mbox{flat system of equations} : \\
                                      ~~~~ c \in e_x \cap A,~~\ x \in X
                                 \end{array} \right\} \subseteq SET $ \\ \\

and clearly, $V [ A ]$ is always a proper class, see 2) in Note 2.1. below. \\

\hfill $\Box$ \\

There are {\it two remarkable} facts about the concept of {\it flat systems of equations} given in the above Definition 2.1., namely

\begin{itemize}

\item the Anti-Foundation Axiom, (AFA), upon which the whole of Self-Referential Mathematics in [2] rests, has a most simple formulation in terms of
{\it flat systems of equations}, as seen next,

\item the {\it flat systems of equations} do in fact give {\it all} the additional new sets in $SET \setminus SET_0$, that is those which due to their
self-referential definitions, cannot be obtained by the usual ZFC Set Theory, see the equivalence Theorem 2.1. below, see also 1) in Note 2.1. below.

\end{itemize}

{~} \\

$\begin{array}{|l}

\mbox{ANTI-FOUNDATION AXIOM (AFA)} \\ \\

{}~~~ \forall~~ {\cal E} ~\mbox{flat system of equations} ~:~ \exists \,!~~ s~~ \mbox{solution} \end{array} $ \\ \\

Here, for the sake of further clarity, let us recall that in ZFC we have \\ \\

$\begin{array}{|l}

\mbox{AXIOM OF FOUNDATION ( FA )} \\ \\

{}~~~ \forall~~ a \in SET ~:~ < a, \in \, > ~\mbox{well-founded} \end{array} $ \\ \\

where the concept of a {\it well-founded} relational structure is defined as follows. \\

A relational structure $ < S, R >$ is called {\it well-founded}, if and only if it has no infinite descending sequence \\

$ \ldots a_n \, R \, a_{n - 1} \, R \ldots R \, a_2 \, R \, a_1 \, R \, a_0$ \\

with $a_0, a_1, a_2, \ldots , a_{n - 1}, a_n, \ldots \in S$. \\

For a relational structure $ < S, R >$, we denote \\

$ < S, R >_{wf} ~=~ \{ a \in S ~|~ \mbox{no infinite descending sequence in}~ R ~\mbox{starting with}~ a \} $ \\

Clearly \\

$ < S, R > well-founded ~~\Longleftrightarrow~~ < S, R >_{wf} ~=~ S $ \\

For simplicity, we denote \\

$R_{wf} \, = \, < S, R>_{wf}$ \\

{\bf Remark 2.1.} \\

In view of the above (AFA) axiom, the question arises :

\begin{itemize}

\item Which flat systems of equations have solutions under the (FA) axiom ?

\end{itemize}

The answer is obtained based on the following concept. A flat system of equations ${\cal E} \, = \, < X, A, e >$ is called {\it well-founded}, if and only
if the relation $<$ defined on $X$ by \\

$~~~ x < y ~~~\Longleftrightarrow~~~ y \in e_x $ \\

is well-founded. And then we have, see [1,2] \\

{\bf Mostowski Collapsing Lemma} \\

In ZFC$^-$ well-founded flat systems of equations have unique solutions. \\

{\bf Corollary 2.1.} \\

In ZFC$^-$ we have the equivalence \\

$(FA) ~~\Longleftrightarrow~~ \left ( ~ \begin{array}{l}
                                             \mbox{ {\it only} well-founded flat systems} \\
                                             ~~\mbox{of equations have solutions}
                                        \end{array} ~ \right )$ \\ \\

{\bf Note 2.1.} \\

1) Flat systems of equations can trivially recover all sets $E \in SET$. Indeed, let $A = E$. Further, let $X = \{ x \}$, with $x \in {\cal U}$, such that
$x \notin E$, which is possible since ${\cal U} \nsubseteqq E$, given the fact that ${\cal U}$ is a proper class. Then $X \subseteq {\cal U}$ and $X
\cap A = \phi$, hence \\

$ e_x = E$ \\

is obviously a flat system of equations, thus according to (AFA), it has a unique solution $s$. Now in view of (2.6),
(2.7), we have $b_x = e_x \cap X = \phi,~ c_x = e_x \cap A = E$, and then (2.9) gives \\

$ s_x = \{ s_y ~|~ y \in e_x \cap X \} \cup ( e_x \cap A ) = E$ \\

In this way, the set $E$ was obtained as the unique solution \\

$ s : X = \{ x \} \ni x \longmapsto s_x = E$ \\

of the above flat system of equations. \\

2) In view of the above example, each $x \in {\cal U}$ leads to a flat system of equations with the respective unique solution $s_x$. And clearly,
if $x, \, x\,' \in {\cal U},~ x \neq x\,'$, then $s_x \neq s_{x\,'}$. As for ${\cal U}$, it is a proper class, therefore, so is $\{ s_x ~|~ x \in
{\cal U} \}$. \\

{\bf Examples 2.1.} \\

Let us illustrate the above in the case of the equations (2.1) - (2.3). \\

1) For (2.1), we can take \\

$ X = \{ x \} \subseteq {\cal U},~~ A = \{ \phi \},~~ e_x = \{ x \} \in {\cal P} ( X \bigcup A )$ \\

therefore, it is a flat system of equations, made up of a single equation. As for its unique solution $s_x \in SET$, we
shall see the details in 4) and 5) below. \\

2) For (2.2), we can take \\

$ X = \{ x \} \subseteq {\cal U},~~ A = \{ \phi \},~~ e_x = x $ \\

thus \\

$e_x \in X \bigcup A,~~ e_x \nsubseteq X \bigcup A,~~ e_x \notin {\cal P} ( X \bigcup A )$ \\

therefore, it is {\it not} a flat system of equations. \\

Also, with (2.2), we can immediately note why the condition \\

(2.10)~~~ $ X \subseteq {\cal U} $ \\

was requested in Definition 2.1. Indeed, without that condition, equation (2.2) is satisfied by {\it all} sets $a \in SET$, thus (2.2) does {\it not}
have a unique solution in $SET$. \\

3) With the equation (2.3), we can take \\

$X = \{ x \} \bigcup x \subseteq {\cal U},~~ A = \phi,~~ e_x = {\cal P} ( x )
                                \nsubseteqq X \bigcup A,~~ e_x  \notin {\cal P} ( X \bigcup A )$ \\

which, however, does not turn (2.3) into a flat system of equations. Also, as seen in 3) in Proposition 2.3. below, equation (2.3) does {\it not} have
any solution in SET. \\

4) In ZFA, the equation \\

(2.11)~~~ $x = \{x \}$ \\

has a {\it unique solution} $\Omega \in SET$. Indeed, as note at 1) above, if we take \\

$X = \{ x \} \subseteq {\cal U},~~ A = \phi,~~ e_x = \{ x \}$ \\

then (2.11) is a flat system of equations, thus in view of (AFA), it has a unique solution $s_x \in SET$, and according to (2.9), we have \\

(2.12)~~~ $s_x = \{ s_y ~|~ y \in b_x \} \cup c_x$ \\

However, (2.6) gives $b_x = e_x \cap X = \{ x \} \cap X = X = \{ x \}$, while (2.7) implies $c_x = e_x \cap A = \phi$.
Thus (2.12) becomes \\

$s_x = \{ s_x \}$ \\

5) The above unique solution $\Omega \in SET$ obviously has the property \\

$\Omega = \{ \Omega \} = \{ \{ \Omega \} \} = \{ \{ \{ \Omega \} \} \} = \ldots$ \\

however, this need {\it not} mean that the bracket pairs $\{~\}$ could be infinitely many, namely, that we could have for instance \\

$\Omega = \ldots \{ \{ \{ \Omega \} \} \} \ldots$ \\

let alone that the bracket pairs $\{~\}$ could reach to transfinite ordinals, or go through all the ordinals, see Remark 2.3. below. \\

6) In ZFA, there is a {\it unique} set \\

(2.13)~~~ $\{ 0, \{ 1, \{ 2, \{ 3, \ldots \} \} \} \} \in SET$ \\

Indeed, we consider the flat system of equations \\

$x_0 = \{ 0, x_1 \}$ \\
$x_1 = \{ 1, x_2 \}$ \\
$x_2 = \{ 2, x_3 \}$ \\
$x_3 = \{ 3, x_4 \}$ \\
\vdots \\

where $X = \{ x_0, x_1, x_2, x_3, \ldots \} \subseteq {\cal U},~~ A = \{ 0, 1, 2, 3, \ldots \}$, while $e_{x_n} = \{ n, x_{n + 1}\}$, with $n \geq 0$.
Then (AFA) gives a unique solution $s$, and in view of (2.6), (2.7), (2.9), we obtain the relations \\

$s_{x_0} = \{ s_x ~|~ x \in b_{x_0} \} \cup c_{x_0} = \{ s_x ~|~ x \in e_{x_0} \cap X \} \cup ( e_{x_0} \cap A ) = \{ 0, s_{x_1} \}$ \\
$s_{x_1} = \{ s_x ~|~ x \in b_{x_1} \} \cup c_{x_1} = \{ s_x ~|~ x \in e_{x_1} \cap X \} \cup ( e_{x_1} \cap A ) = \{ 1, s_{x_2} \}$ \\
$s_{x_2} = \{ s_x ~|~ x \in b_{x_2} \} \cup c_{x_2} = \{ s_x ~|~ x \in e_{x_2} \cap X \} \cup ( e_{x_2} \cap A ) = \{ 2, s_{x_3} \}$ \\
\vdots \\

thus \\

$s_{x_0} = \{ 0, s_{x_1} \} = \{ 0,  \{ 1, s_{x_2} \} \} =  \{ 0,  \{ 1, \{ 2, s_{x_3} \} \} \} = \ldots$ \\

7) Let us consider that flat system of equations {\it without} atoms, that is, with $A = \phi$, namely \\

$x_0 = \{ y_0, x_1 \}~~~~ y_0 = \phi$ \\
$x_1 = \{ y_1, x_2 \}~~~~ y_1 = \{ y_0 \}$ \\
$x_2 = \{ y_2, x_3 \}~~~~ y_2 = \{ y_0, y_1 \}$ \\
$x_3 = \{ y_3, x_4 \}~~~~ y_3 = \{ y_0, y_1, y_2 \}$ \\
\vdots \\

where $X = \{ x_0, y_0, x_1, y_1, x_2, y_2, \ldots \} \subseteq {\cal U},~~ A = \phi$, while $e_{x_n} = \{ y_n, x_{n+ 1} \},~ e_{y_{n + 1}} = \{ y_0,
\ldots , y_n \}$, for $n \ge 0$, and $e_{y_0} = \phi$. In this case, (2.6), (2.7) give \\

$b_{x_n} = e_{x_n} \cap X = \{ y_n, x_{n + 1} \},~~~ n \geq 0$ \\

$b_{y_0} = e_{y_0} \cap X = \phi$ \\
$b_{y_{n + 1}} = e_{y_{n + 1}} \cap X = \{ y_0, \ldots , y_n \},~~~ n \geq 0$ \\

$c_{x_n} = c_{y_n} = \phi,~~~ n \geq 0$ \\

therefore, in view of (2.9), the unique solution $s$ given by (AFA), is such that \\

$s_{x_n} = \{ s_x ~|~ x \in b_{x_n} \} \cup c_{x_n} = \{ s_{y_n}, s_{x_{n + 1}} \},~~~ n \geq 0$ \\

$s_{y_0} = \{ s_y ~|~ y \in b_{y_0} \} \cup c_{y_0} = \phi$ \\
$s_{y_{n + 1}} = \{ s_y ~|~ y \in b_{y_{n + 1}} \} \cup c_{y_{n + 1}} = \{ s_{y_0}, \ldots , s_{y_n} \},~~~ n \geq 0$ \\

In particular \\

$s_{y_0} = \phi$ \\
$s_{y_1} = \{ s_{y_0} \} = \{ \phi \}$ \\
$s_{y_2} = \{ s_{y_0}, s_{y_1} \} = \{ \phi, \{ \phi \} \}$ \\
$s_{y_3} = \{ s_{y_0}, s_{y_1}, s_{y_2} \} = \{ \phi, \{ \phi \}, \{ \phi, \{ \phi \} \} \}$ \\
\vdots \\

which means that, for $n \geq 0$, we obtain that \\

$s_{y_n} = n$~ in the von Neumann representation \\

{\bf Proposition 2.1.} \\

Let ${\cal E} =\, < X, A, e >$\, be a flat system of equations. If $A \subseteq {\cal U}$, then set $solution-set
({\cal E})$ is {\it transitive}, namely \\

$ b,~c \in SET,~ c \in b \in solution-set ( {\cal E} ) ~~~\Longrightarrow~~~ c \in solution-set ( {\cal E} )$ \\

{\bf Proof}. \\

We recall that \\

$solution-set ({\cal E}) = \{ s_y ~|~ y \in e_x \cap X,~ x \in X \} \, \bigcup \, \bigcup_{x \in X} ( e_x \cap A )$ \\

Let $b \in solution-set ( {\cal E} )$, then \\

either $b = s_y$, for some $y \in e_x \cap X$, with suitable $x \in X$, \\

or $b \in e_x \cap A$, for some $x \in X$. \\

In the first case, if $c \in b = s_y = \{ s_z ~|~ z \in e_y \cap X \} \cup ( e_y \cap A )$, then \\

either $c = s_z$, thus $c \in solution-set ({\cal E})$, \\

or $c \in e_y \cap A $, thus again  $c \in solution-set ({\cal E})$. \\

In the second case, if $b \in e_x \cap A$, then $b \in {\cal U}$, thus there cannot be $c \in b$. \\

{\bf Proposition 2.2.} \\

If $A \subseteq {\cal U}$, then \\

$~~~ V [ A ] \subseteq V_{afa} [ A ] $ \\

where the operation $V_{afa}$ is defined in (3.5) in the next section. \\

{\bf Note :} Here we make an advance use of some notations and results in section 3 below. However, placing Proposition 2.2. here in section 2 helps in
the better understanding of the concept of flat system of equations, as well as of its fundamental importance seen in the equivalence Theorem 2.1.
below. \\

{\bf Proof}. \\

Let $c \in  V [ A ]$, then for some ${\cal E} =\, < X, A, e >$ flat system of equations we have either \\

$ c = s_y$, for some $y \in e_x \cap X,~ x \in X$ \\

or \\

$ c \in e_x \cap A$, for some $x \in X$ \\

In the first case we also have \\

$ c \subseteq Z = solution-set ({\cal E}) =$ \\

$~~~~~~~~~~ = \{ s_y ~|~ y \in e_x \cap X,~ x \in X \} \, \bigcup \, \bigcup_{x \in X} ( e_x \cap A )$ \\

since in view of (2.9) \\

$ c = s_y = \{ s_z ~|~ z \in e_y \cap X \} \cup ( e_y \cap A )$ \\

But in view of Proposition 2.1., the set $Z = solution-set ({\cal E})$ is transitive. Therefore, for $x \in X$, we have,
see (3.3) below \\

$ s_x \subseteq TC ( s_x ) \subseteq Z$ \\

which gives \\

$ TC ( s_x ) \cap {\cal U} \subseteq Z \cap {\cal U} \subseteq A$ \\

Indeed, if $b \in Z \cap {\cal U}$, then in particular \\

either $b = s_y$, for some $y \in e_x \cap X$, with suitable $x \in X$, \\

or $b \in e_x \cap A$, for some $x \in X$. \\

In the first case, $b \notin {\cal U}$, since $s_y \in SET$. \\

In the second case obviously $b \in A$. \\

In conclusion, in view of (3.5), we have $c \in V_{afa} [ A ]$. \\

{\bf Remark 2.1.} \\

The following is, of course, a fundamental question :

\begin{itemize}

\item How many sets $a \in SET$ can be obtained as solutions of flat systems of equations ?

\end{itemize}

In 1) in Note 2.1. above, we have seen that flat systems of equations can trivially recover as solutions all sets in $SET$. A more precise and rather
natural answer, and as such, best possible answer will be given in Theorem 2.2. below. \\
Needless to say, this answer highlights the importance of flat systems of equations. However, in various contexts, other two concepts of systems of
equations will prove to be useful, concepts given in Definitions 2.2. and 2.3. below.

\hfill $\Box$ \\

Now, the second kind of systems of equations aims to eliminate the above restriction $X \subseteq {\cal U}$ in (2.10) on
the flat systems of equations. And as we shall see in Theorem 2.1. below, this is in fact possible, in spite of the above problem with
lack of uniqueness of solutions, provided that the STRONG AXIOM OF PLENITUDE is accepted. \\

{\bf Definition 2.2.} \\

A structure ${\cal E} ~=~ < X, A, e >$ is called a {\it generalized flat system of equations}, if and only if \\

(2.14)~~~ $ X,~ A \in SET,~~~ X \bigcap A = \phi $ \\

with $X$ the set of $indeterminates$ and $A$ the set of $atoms$, while \\

(2.15)~~~ $ e : X \longrightarrow {\cal P} ( X \bigcup A ) $

\hfill $\Box$ \\

There is a close connection between the solutions of flat, and on the other hand, generalized flat systems of equations, provided that the following
axiom holds \\ \\

$\begin{array}{|l}

\mbox{STRONG AXIOM OF PLENITUDE} \\ \\

\mbox{There is an operation}~ new ( a, b ), ~\mbox{such that} \\ \\

{}~~~ 1) ~~~\forall~~ a \in SET,~ b \subset {\cal U} ~:~ new ( a, b ) \in {\cal U} \setminus a \\ \\

{}~~~ 2) ~~~\forall~~ a, a\,' \in SET,~ a \neq a\,',~ b \subset {\cal U} ~:~ new ( a, b ) \neq new ( a\,', b ) \end{array} $ \\ \\

{\bf Theorem 2.1.} \\

Assuming the STRONG AXIOM OF PLENITUDE, every generalized flat system of equations ${\cal E} ~=~ < X, A, e >$  has a
unique solution $s$. Furthermore, there exists an associated flat system of equations ${\cal E}\,' ~=~ < Y, A, e\,' >$,
such that \\

$solution-set ( {\cal E} ) = solution-set ( {\cal E}\,' )$ \\

{\bf Proof}. \\

We have to replace $X$ by a set $Y \subset {\cal U}$, such that $Y \cap A =\phi$. Thus we take \\

$Y = \{ y_x ~|~ y_x = new ( x, A ),~~~ x \in X \}$ \\

then \\

$Y \subset {\cal U},~~ Y \cap A = \phi$ \\

Let now \\

$e\,'_{y_x} = \{ y_z ~|~ z \in e_x \cap X \} \bigcup \,( e_x \cap A )$ \\

Clearly, ${\cal E}\,' ~=~ < Y, A, e\,' >$ is a flat system of equations, and thus, it has a unique solutio $s\,'$. \\

Now we get the solution $s$ of ${\cal E} ~=~ < X, A, e >$ given by \\

$s_x = s\,'_{y_x},~~~ x \in X$ \\

The uniqueness of $s$ follows from the fact that every solution $t$ of ${\cal E} ~=~ < X, A, e >$ gives a solution $t\,'$
of ${\cal E}\,' ~=~ < Y, A, e\,' >$. And we must have $t\,' = s\,'$, thus it follows that $t = s$.

\hfill $\Box$ \\

{\bf Example 2.2.} \\

For every set $a \in SET$, we associate the {\it canonical generalized flat system of equations} ${\cal E}_a = \,
< X_a, A_a, e_a >$, with, see (3.1) \\

$ A_a = TC ( a ) \bigcap {\cal U}$ \\

$ X_a = TC ( \{ a \} ) \setminus A_a $ \\

and \\

$ X_a \ni x \longmapsto ( e_a )_x = x \subseteq X_a \bigcup A_a$ \\

where we have to prove the inclusion in the last relation. \\

We note in this regard that, see (3.2$^*$), (3.3$^*$) and 1) in Examples 3.1.  \\

$ A_a = \{ x \in {\cal U} ~|~ x \in a_n \in \ldots \in a_1 \in a_0 = a,~~ n \geq 0 \}$ \\

$ X_a = (~ \{ a \} \, \bigcup  \, TC ( a )  ~) \setminus (~ TC ( a ) \bigcap {\cal U} ~) =$ \\

$~~~~~~~~~~~~ = (~ \{ a \} \, \bigcup \, \{ x ~|~ x \in a_n \in \ldots \in a_1 \in a_0 = a,~~ n \geq 0 \} ~) \setminus$ \\

$~~~~~~~~~~~~~~~~~~ \setminus \{ y \in {\cal U} ~|~ y \in b_m \in \ldots \in b_1 \in b_0 = a,~~ m \geq 0 \} =$ \\

$~~~~~~~~~~~~ = \{ a \} \, \bigcup \, \{ x \in SET ~|~ x \in a_n \in \ldots \in a_1 \in a_0 = a,~~ n \geq 0 \}$ \\

thus \\

$ x \in X_a ~\Longrightarrow~ x \in SET$ \\

Furthermore \\

$ x \in X_a ~\Longrightarrow~ x \subseteq X_a \bigcup A_a$ \\

Indeed, given $x \in X_a$, then there are the two cases \\

(i.1)$~~~ x = a$ \\

(i.2)$~~~  x \in a_n \in \ldots \in a_1 \in a_0 = a$, for some $n \geq 0$ \\

Let now $y \in x$. In case (i.1), we have $y \in a$, thus $y \in X_a$, provided that $y \notin {\cal U}$. Otherwise obviously $y \in A_a$. In case (i.2),
clearly $y \in ( X_a \cup A_a )$. \\

Assuming now the STRONG AXIOM OF PLENITUDE, we have in view of Theorem 2.1., a unique solution $s_a$ of ${\cal E}_a$. And in fact, we have \\

$ s_a : X_a \ni x \longmapsto ( s_a )_x = x \in SET$ \\

Indeed, (2.9) gives for $x \in X_a$ \\

$ ( s_a )_x = \{ ( s_a )_y ~|~ y \in ( e_a )_x \cap X_a \} \, \cup \, ( ( e_a )_x \cap A_a ) \in SET $ \\

while \\

$ ( e_a )_x \cap X_a = x \cap X_a,~~~ ( e_a )_x \cap A_a = x \cap A_a$ \\

thus \\

$ ( s_a )_x = \{ ( s_a )_y ~|~ y \in x \cap X_a \} \, \cup \, ( x \cap A_a ) \in SET $ \\

which in our case becomes the identity \\

$ x = \{ y ~|~ y \in x \cap X_a \} \, \cup \, ( x \cap A_a ) \in SET $ \\

{\bf Theorem 2.2. Equivalence} \\

Assuming now the STRONG AXIOM OF PLENITUDE, we have for $A \subseteq {\cal U}$ the equivalence between sets in $SET$ which have support in $A$, and
sets in $SET$ which are solutions of flat systems of equations with atoms in $A$, namely \\

$~~~~~~ V_{afa} [ A ] =  V [ A ] $ \\

{\bf Proof}. \\

In view of Proposition 2.2., we have \\

$ V [ A ] \subseteq V_{afa} [ A ]$ \\

Let now $a \in V_{afa} [ A ]$. Then by definition, see (3.5) below \\

$ TC ( a ) \cap {\cal U} \subseteq A$ \\

Now, in view of Example 2.2., we consider the unique solution $s_a$ of the canonical generalized flat system of equations ${\cal E}_a = \,
< X_a, A_a, e_a >$, which is $s_a = id_{X_a}$. Thus, recalling that $a \in X_a$ and $A_a =  TC ( a ) \cap {\cal U} \subseteq A$, we obtain \\

$ s_a = a$ \\

which gives \\

$ a \in solution-set ( {\cal E}_a )$ \\

On the other hand, in view of Theorem 2.1., there is a flat system of equations ${\cal E} = \, < X, A_a, e >$ with the same atoms $A_a$, such that \\

$ solution-set ( {\cal E}_a ) = solution-set ( {\cal E} )$ \\

hence $a \in solution-set ( {\cal E} ) \subseteq V [ A ]$

\hfill $\Box$ \\

Finally, the third kind of systems of equations allows considerably more general right hand terms $e_x$ in (2.5), (2.15), although it has to accept
harder restrictions on $X$ and $A$, than in Definition 2.1., that is, in (2.10). Namely \\

{\bf Definition 2.3.} \\

A structure ${\cal E} ~=~ < X, A, e >$ is called a {\it generalized system of equations}, if and only if \\

(2.16)~~~ $ X,~ A \in SET,~~~ X,~ A \subseteq {\cal U},~~~ X \bigcap A = \phi $ \\

with $X$ the set of $indeterminates$ and $A$ the set of $atoms$, while \\

(2.17)~~~ $ e : X \longrightarrow V_{afa} ( X \bigcup A ) $ \\

where the operation $V_{afa}$ is defined in (3.5) in the next section. \\

{\bf Remark 2.2.} \\

As we shall see in section 3, the range $V_{afa} ( X \bigcup A )$ of the mappings $e$ in (2.17) is considerably {\it
larger} than ${\cal P} ( X \bigcup A )$, which is the range of the corresponding mappings in (2.5) and (2.15). Therefore,
the generalized systems of equations defined above contain as a rather {\it small} particular case the flat and the
generalized flat systems of equations. \\

{\bf Examples 2.3.} \\

1) For a given $a \in SET \bigcup {\cal U}$, let us consider the equation \\

(2.18)~~~ $x = x \longrightarrow a$ \\

hence it is {\it not} a flat or generalized flat system of equations, if we take \\

$X = \{ x \},~~ A =\{ a \},~~ e_x =  x \longrightarrow a \subseteq x \times a,~~ e_x \nsubseteq X \bigcup A$ \\

$~~~~~~~~~~ e_x \notin {\cal P} ( X \bigcup A )$ \\

although (2.18) {\it has} solution in SET, see 5) in Proposition 2.3. below. \\

2) For given $p, q \in SET \bigcup {\cal U}$, the equation \\

(2.19)~~~ $x = \{ \{ x, q \}, p \}$ \\

is {\it not} a flat or generalized flat system of equations, if considered with \\

$X = \{ x \},~~ A = \{ p, q \},~~ e_x = \{ \{ x, q \}, p \} \nsubseteq X \bigcup A,~~ e_x \notin {\cal P} ( X \bigcup A )$ \\

although it can be written as a flat system of equations, provided that $x \in {\cal U}$, namely \\

{}~~~~~ $x = \{ y, p \}$ \\

{}~~~~~ $y = \{ x, q \}$ \\

hence \\

$X = \{ x, y \},~~ A = \{ p, q \},~~ e_x =  \{ y, p \}, e_y = \{ x, q \} \subseteq X \bigcup A$, \\

$e_x, e_y \in {\cal P} ( X \bigcup A )$ \\

{\bf Proposition 2.2.} \\

Within ZFC we have \\

1) $~~~ \forall~~ a \in SET ~:~ a \notin a $ \\

2) $~~~ \neg \exists~~ a_1, \ldots , a_n \in SET ~:~ a_1 \in \ldots a_n \in a_1 $ \\

3) $~~~ \neg \exists~~ a, b \in SET ~:~ a \in TC ( b ) \in a $ \\

4) $~~~ \neg \exists~~ a, b \in SET ~:~ a \in TC ( b ),~~ b \in TC ( a ) $ \\

5) $~~~ \forall~~ a, b, c \in SET ~:~ c \,=\, < a, b > ~\Longrightarrow~ c \neq a,~ c \neq b,~ c \notin a,~ c \notin b $ \\

6) $~~~ \forall~~ A, X \in SET ~:~ X \neq \phi ~\Longrightarrow~ X \neq A \times X $ \\

7) $~~~ \forall~~ X \in SET ~:~ X = X \times X ~\Longrightarrow~ X = \phi $ \\

8) $~~~ \neg \exists~~ \mbox{function}~~ f = A \longrightarrow B,~ A, B \in SET ~:~ f \in dom ( f ) $ \\

9) $~~~ \neg \exists~~ \mbox{functions}~~ f_1 : A_1 \longrightarrow A_2, \ldots , f_n : A_n \longrightarrow A_1,~\\ \\
        \hspace*{2cm} A_1, \ldots , A_n \in SET,~ a_1 \in A_1 ~:~ \\ \\
        \hspace*{2.5cm} f_n ( \ldots f_1( a_1 ) \ldots ) = f_1 $ \\

10) $~~~ \forall~~ A, X \in SET ~:~ X \neq X \longrightarrow A $ \\

It is important to note that, as seen next, even in ZFC$^-$, that is, without FA, one can obtain {\it impossibilities} of
self-reference. \\

{\bf Proposition 2.3.} \\

Within ZFC$^-$ we have \\

1) $~~~ \forall~ \mbox{function}~ F : A \longrightarrow B,~ A, B \in SET ~:~ \\ \\
        \hspace*{2cm} \{ x \in dom ( F ) ~|~ x \notin F ( x ) \} \notin rng ( F ) $ \\

2) $~~~ \forall~ \mbox{function}~ F : A \longrightarrow {\cal P} ( A ),~ A\in SET ~:~ \\ \\
        \hspace*{2cm} \{ x \in A ~|~ x \in R_{wf} \} \notin rng ( F ) $ \\

{~~~~~~} where $R = \{ < x, y > ~|~ x \in F ( y ) \}$ \\

3) $~~~ \forall~~ X \in SET ~:~ X \neq {\cal P} ( X ) $ \\

4) $~~~ \forall~~ X \in SET ~:~ X = X \longrightarrow \phi ~~\Longrightarrow~~ X = \{ \phi \} $ \\

5) $~~~ \forall~~ A, X \in SET ~:~ X = X \longrightarrow A ~~\Longrightarrow~~ \\ \\
        \hspace*{2.5cm} \Longrightarrow~~ A = \{ a \},~ X = \{ f \},~ f ( f ) = a $ \\

6) $~~~ \forall~~ X \in SET ~:~ X = X \longrightarrow X ~~\Longrightarrow~~ X = \{ x \},~ x = \{ < x, x > \} $ \\

The proofs of the above two Propositions 2.2. and 2.3. are rather simple and immediate, and can be found at [1, pp. 25-27] \\

{\bf Remark 2.3.} \\

Related to 5) and 6) in Examples 2.1. above, let us note the following two kind of situations encountered so far with sets which have {\it infinitely
many brackets}, namely : \\

(2.20)~~~ $\Omega \stackrel{?}{=} \ldots \{ \{ \{ \Omega \} \} \} \ldots$ \\

(2.21)~~~ $\{ 0, \{ 1, \{ 2, \{ 3, \ldots \} \} \} \} \in SET$ \\

The second one was, in 6) in Examples 2.1. above, proved to exist uniquely, and be well defined in ZFA, while the first one will be considered in more
detail in 2) below. \\

1) Related to (2.21), we note the following immediate generalization. Let $\alpha$ be any infinite ordinal number and let us take \\

$X = \{ x_\beta ~|~ \beta < \alpha \} \subseteq {\cal U},~~ A = \{ \beta ~|~ \beta < \alpha \}$ \\

while \\

$e_{x_\beta} = \{ \beta, x_{\beta + 1} \},~~~ \beta < \alpha $ \\

Then obviously, we obtain a flat system of equations, therefore (AFA) gives a unique solution $s$ which, in view of (2.9), has the property \\

$~~~~~~~~~~ s_{x_\beta} = \{ s_x ~|~ x \in b_{x_\beta} \} \bigcup c_{x_\beta} \in SET,~~~ \beta < \alpha$ \\

where according to (2.6), (2.7), we have \\

$b_{x_\beta} = e_{x_\beta} \cap X = \{ x_{\beta + 1} \},~~~ c_{x_\beta} = e_{x_\beta} \cap A = \phi,~~~ \beta < \alpha$ \\

hence \\

$~~~~~~~~~~ s_{x_\beta} = \{ s_{x_{\beta + 1}} \} \in SET,~~~ \beta < \alpha$ \\

which gives \\

(2.22)~~~ $s_{x_0} = \{ 0, s_{x_1} \} = \{ 0,  \{ 1, s_{x_2} \} \} =  \{ 0,  \{ 1, \{ 2, s_{x_3} \} \} \} = \ldots$ \\

where the pairs of brackets $\{~\}$ occur once for each $\beta < \alpha$. \\

Thus (2.21) is the particular case of the above $s_{x_0}$ in (2.22) corresponding to $\alpha = \omega$ which is the first infinite ordinal. In the general
case of an infinite ordinal $\alpha$, the above $s_{x_0}$ in (2.22) gives instead of (2.21) the set \\

(2.23)~~~ $\{ 0, \{ 1, \{ 2, \{ 3, \ldots \{ \beta, \ldots \} \ldots \} \} \} \} \in SET$ \\

which contains all $\beta < \alpha$. \\

2) Let us return to (2.20) and consider it as a particular case of the following general operation : given $a \in SET \bigcup {\cal U}$,
define the set \\

(2.24)~~~ $\ldots \{ \{ \{ a \} \} \} \ldots \in SET$ \\

with a pair of brackets $\{~\}$ for each $n < \infty$. For that purpose, let use the following notation \\

$\underset{0}\{ a \underset{0}\} ~=~ a,~~~ \underset{1}\{ a \underset{1}{\}} ~=~ \{ a \},~~~
                    \underset{2}\{ a \underset{2}\} ~=~ \{ \{ a \} \}, \ldots $ \\

Thus the problem is :

\begin{quote}

How to define in $SET$

\end{quote}

(2.25)~~~ $\underset{\omega}\{ a \underset{\omega}\} ~ \in SET$

\begin{quote}

where $\omega$ denotes the first infinite ordinal number.

\end{quote}

Of course, one would want to define (2.25) as a certain kind of "limit" of the sequence of sets $\underset{n}\{ a \underset{n}\} \, \in SET,~~ n
\geq 0$. \\

One way to do that for an arbitrary set $a \in SET$ is as follows. Let us denote \\

$ \overset{0}\{ a \overset{0}\} ~=~ a$ \\

$ \overset{1}\{ a \overset{1}{\}} ~=~ \overset{0}\{ a \overset{0}\} \cup \{ a \} = a \cup \{ a \}$ \\

$ \overset{2}\{ a \overset{2}\} ~=~ \overset{1}\{ a \overset{1}{\}} \cup \{ a, \{ a \} \} = a \cup \{ a \} \cup \{ a, \{ a \} \}$ \\

$ \overset{3}\{ a \overset{3}\} ~=~ \overset{2}\{ a \overset{2}\} \cup \{ a, \{ a, \{ a \} \} \} =
                     a \cup \{ a \} \cup \{ a, \{ a \} \} \cup \{ a, \{ a, \{ a \} \} \}$ \\

$ \overset{4}\{ a \overset{4}\} ~=~ \overset{3}\{ a \overset{3}\} \cup \{ a, \{ a, \{ a, \{ a \} \} \} \} = $ \\

$~~~~~~~~~~~~ = a \cup \{ a \} \cup \{ a, \{ a \} \} \cup \{ a, \{ a, \{ a \} \} \} \cup \{ a, \{ a, \{ a, \{ a \} \} \} \}$ \\
\vdots \\

It follows that \\

(2.26)~~~ $\overset{1}\{ a \overset{1}{\}} \subseteq \overset{2}\{ a \overset{2}\} \subseteq \overset{3}\{ a \overset{3}\}
                              \subseteq \overset{4}\{ a \overset{4}\} \ldots $ \\

Thus one can define \\

(2.27)~~~ $\overset{\omega}\{ a \overset{\omega}\} ~=~ \bigcup_{n < \omega} \overset{n}\{ a \overset{n}\} \in SET$ \\

Clearly, that procedure can be extended to all ordinal numbers $\alpha$. \\

Thus the above problem (2.25) got solved in general, although not along its initial formulation. \\

On the other hand, in the particular case when $a = \Omega \in SET$, then in view of the fact that \\

$ \Omega \,=\, \underset{n}\{ a \underset{n}\} \, \in SET,\,~~ n < \omega$ \\

one may come up with a definition of (2.25) considered in its initial formulation, and which hence is simpler than the one given in (2.27), namely \\

(2.28)~~~ $\underset{\omega}\{ a \underset{\omega}\} \,=\, \Omega \in SET$ \\

And again, one may extend that definition to all ordinal numbers $\alpha$, by \\

(2.29)~~~ $\underset{\alpha}\{ a \underset{\alpha}\} \,=\, \Omega \in SET$ \\

3) The obvious difference between (2.20) and (2.21) is that in the second, there is an {\it outer} pair of brackets $\{~\}$, while in the first there is
none. And such an outer pair of brackets does indeed define a set in $SET$, or for that matter, even in $SET_0$, provided that what is within that outer
pair of brackets makes sense in the respective version of Set Theory. And clearly, for (2.21) such is the case within $SET$, as seen in 6) in Examples
2.1. above. \\

One can also note that the generalization of (2.21) in (2.23) to arbitrary ordinals $\alpha$ always has an outer pair of brackets $\{~\}$. On the
other hand, in the generalization (2.29) of (2.20), there is an outer pair of brackets $\{~\}$, only if $\alpha$ is not a limit ordinal. \\

4) The flat system of equations in 7) in Examples 2.1., can obviously be generalized to arbitrary ordinal numbers $\alpha$, in a way similar to the
generalization in 1) above of 6) in Examples 2.1. \\ \\

{\bf 3. Three Basic Operations} \\

In order to pursue the theory, the following three operations, seldom if at all encountered in usual Set Theory, although
quite elementary as such, will be needed. \\

We start with the definition of a fundamental concept. \\

{\bf Definition 3.1.} \\

A set $a \in SET$ is called {\it transitive}, if and only if \\

$~~~~~~~~~~~~ b \in a ~~\Longrightarrow~~ b \subseteq a$ \\

or equivalently \\

$~~~~~~~~~~~~ c \in b \in a ~~\Longrightarrow~~ c \in a$  \\

Clearly, usual sets in mathematics are {\it not} transitive. For instance, given a set ${\cal X}$ of open subsets in a topological space, then the
transitivity of ${\cal X}$ would imply that for every open subset $E \in {\cal X}$, we must also have $E \subseteq {\cal X}$. In other words, ${\cal X}$
must also contain as elements all the points $x \in E$, for every $E \in {\cal X}$. \\

And now, the first basic operation. \\

{\bf Definition 3.2.} \\

Given a set $a \in SET$, its {\it transitive closure} is by definition the smallest transitive set which contains it, and which is denoted by
$TC ( a )$. \\

{\bf Lemma 3.1.} \\

The transitive closure $TC ( a )$ exists for every set $a \in SET$, and it is given by \\

(3.1)~~~ $ TC ( a ) = \bigcup~ \{ a, \bigcup a, \bigcup \bigcup a, \ldots \} \in SET $ \\

Further, for $a \in SET$, we have \\

(3.2)~~~ $TC ( a ) = \{ b ~|~ b \in a \} ~\bigcup~ \{ c \in b \in a \} ~\bigcup~ \{ d \in c \in b \in a \} ~\bigcup$ \\

$~~~~~~~~~~~~~ \bigcup~ \{ e \in d \in c \in b \in a \} ~\bigcup~ \ldots$ \\

Here we used the simplifying notation \\

$\{ c \in b \in a \} = \{~ c ~|~ \exists~~ b \in a ~:~ c \in b ~\}$ \\

$\{ d \in c \in b \in a \} = \{~ d ~|~ \exists~~ b \in a ~:~ \exists~~ c \in b ~:~ d \in c ~\}$ \\

$\vdots$ \\

{\bf Note.} \\

The meaning of $TC ( a )$, for a given set $a \in SET$, is clear from (3.2) which, obviously, can be written in the equivalent form \\

(3.2$^*$)~~~ $TC ( a ) = \{ x ~|~ x \in a_n \in \ldots \in a_2 \in a_1 \in a_0 = a,~~ n \geq 0 \}$ \\

{\bf Proof.} \\

We note that \\

$\bigcup a = \bigcup_{\, b \in a} b = \{ c \in b \in a \}$ \\

thus \\

$\bigcup \bigcup a = \bigcup_{\, b \in \bigcup a} b = \{ c \in b \in \bigcup a \} = \{ d \in c \in b \in a \}$ \\

and so on \ldots \\

Therefore \\

$TC ( a ) = \{ x \in y \in \{ a, \bigcup a, \bigcup \bigcup a, \ldots \} \} =$ \\

$~~~~~~~~ = \{ x \in y = a \} ~\bigcup~ \{ x \in y = \bigcup a \} ~\bigcup~ $ \\

$~~~~~~~~ ~\bigcup~ \{ x \in y = \bigcup \bigcup a \} ~\bigcup~ \ldots = $ \\

$~~~~~~~~ = a ~\bigcup~ (~ \bigcup a ~) ~\bigcup~ (~ \bigcup \bigcup a ~) ~\bigcup~ \ldots = $ \\

$~~~~~~~~ = a ~\bigcup~ (~ \bigcup_{\, b \in a} b ~) ~\bigcup~ (~ \bigcup_{\, c \in \bigcup a} c ~)
                                    ~\bigcup~ (~ \bigcup_{\, d \in \bigcup \bigcup a} d ~) \ldots = $ \\

$~~~~~~~~ = \{ b \in a \} ~\bigcup~ \{ c \in b \in a \} ~\bigcup~ \{ d \in c \in b \in a \} \ldots $ \\

{\bf Examples 3.1.} \\

1) Given $a \in SET$, then $TC ( \{ a \} )$ is the smallest transitive set which has $a \in SET$ as an element, since \\

$TC ( \{ a \} ) = \{ a \} ~\bigcup~ \{ c \in b \in \{ a \} \} ~\bigcup~
                                             \{ d \in c \in b \in \{ a \} \} \ldots = $ \\

$~~~~~~~~ = \{ a \} ~\bigcup~ \{ c \in b = a \} ~\bigcup~ \{ d \in c \in b = a \} \ldots = $ \\

$~~~~~~~~ = \{ a \} ~\bigcup~ \{ b \in a \} ~\bigcup~ \{ c \in b \in a \} \ldots = \{ a \} \, \bigcup \, TC ( a )$ \\

2) $TC ( \phi ) = \phi$ \\

3) If $a \in {\cal U}$, then $TC ( \{ a \} ) = \{ a \}$ \\

Note : if $a \in {\cal U}$, then $TC ( a )$ is {\it not} defined, since $ a \notin SET$ \\

4) If $a \subseteq {\cal U},~ a \in SET$, then $TC ( a ) = a$ \\

5) If $a \in SET$, then \\

$TC ( a ) = \phi ~\Longrightarrow~ \{ b \in a \} = \phi  ~\Longrightarrow~ a = \phi$ \\

thus \\

$TC ( a ) = \phi ~\Longleftrightarrow~ a = \phi$ \\

6) If $a \subseteq {\cal U},~ a \in SET$, then \\

$TC ( a ) = \phi ~\Longrightarrow~ a = \phi \in SET$ \\

$TC ( a ) = \phi ~\Longleftrightarrow~ a = \phi$ \\

7) If $A, B \in SET, A, B \neq \phi, a \in A, b \in B$, then \\

$< a, b >\, = \{ \{ a \}, \{ a, b \} \} \in A \times B \in SET$ \\

and \\

$TC ( < a, b > ) = \, < a, b > ~\bigcup~ \{ d \in c \in \, < a, b > \} ~\bigcup$ \\

$~~~~~~ ~\bigcup~ \{ e \in d \in c \in \, < a, b > \} ~\bigcup~ \ldots =$ \\

$~~~~~~ =  \{ \{ a \}, \{ a, b \} \} ~\bigcup~ \{ d \in c = \{ a \} \bigvee d \in c = \{ a, b \} \} ~\bigcup$ \\

$~~~~~~ ~\bigcup~ \{ e \in d \in c = \{ a \} \bigvee e \in d \in c = \{ a, b \} \} ~\bigcup~ \ldots =$ \\

$~~~~~~ = \{ \{ a \}, \{ a, b \} \} ~\bigcup~ \{ d \in \{ a , b \} \} ~\bigcup~ \{ e \in d \in \{ a, b \} \} ~\bigcup~ \ldots =$ \\

$~~~~~~ = \{ \{ a \}, \{ a, b \} \} ~\bigcup~ \{ a , b \} ~\bigcup~ \{ e \in \{ a, b \} \} ~\bigcup~ \ldots =$ \\

$~~~~~~ = \{ \{ a \}, \{ a, b \} \} ~\bigcup~ TC ( \{ a, b \} ) =  ~< a, b > ~\bigcup~ TC ( \{ a, b \} )$ \\

8) If $A, B \in SET, A, B \neq \phi, R \subseteq A \times B$, then \\

$TC ( R ) = R ~\bigcup~ \{ c \in \, < a, b > \, \in R \} ~\bigcup~ \{ d \in c \in \, < a, b > \, \in R \} ~\bigcup$ \\

$~\bigcup~ \{ e \in d \in c \in \, < a, b > \, \in R \} ~\bigcup~ \ldots =$ \\

$= R \cup \{ c = \{ a \} | < a, b > \, \in R \} \cup~ \{ c = \{ a, b \} | < a, b > \, \in R \} \cup$ \\

$\cup \{ d \in c = \{ a \} | < a, b > \, \in R \} \cup \{ d \in c = \{ a, b \} | < a, b > \, \in R \} \cup$ \\

$\cup \{ e \in d \in c = \{ a \} | < a, b > \, \in R \} \cup$ \\

$\cup \{ e \in d \in c = \{ a, b \} | < a, b > \, \in R \} \cup  \ldots =$ \\

$= R \cup \{ \{ a \} \, | < a, b > \, \in R \} \cup \{ \{ a, b \} \, | < a, b > \, \in R \} \cup$ \\

$\cup \{ a \, | < a, b > \, \in R \} \cup \{ b \, | < a, b > \, \in R \} \cup$ \\

$\cup \{ e \in a \, | < a, b > \, \in R \} \cup \{ e \in b \, | < a, b > \, \in R \} \cup \ldots$ \\

9) The above goes in particular when $R$ is a function $f : A \longrightarrow B$, or when $R = A \times B$, and in the last case we obtain \\

$TC ( A \times B ) = ( A \times B ) \cup \{ \{ a \} \, | \, a \in A \} \cup \{ \{ a, b \} \, | \, a \in A, b \in B \} \cup$ \\

$\cup \{ a \, | \, a \in A \} \cup \{ b \, | \, b \in B \} \cup$ \\

$\cup \{ e \in a \, | \, a \in A \} \cup \{ e \in b \, | \, b \in B \} \cup \ldots =$ \\

$= ( A \times B ) \cup \{ \{ a \} \, | \, a \in A \} \cup \{ \{ a, b \} \, | \, a \in A, b \in B \} \cup TC ( A ) \cup TC ( B )$ \\

10) If $X \in SET$, then \\

$TC ( {\cal P} ( X ) ) = {\cal P} ( X ) ~\bigcup~ \{ c ~|~ c \in b \in {\cal P} ( X ) \} ~\bigcup$ \\

$~~~~~~ ~\bigcup~ \{ d ~|~ d \in c \in b \in {\cal P} ( X ) \} ~\bigcup~ \ldots =$ \\

$~~~~~~ = {\cal P} ( X ) ~\bigcup~ \{ c ~|~ c \in Y \subseteq X \} ~\bigcup~ \{ d ~|~ d \in c \in Y \subseteq X \} ~\bigcup~ \ldots =$ \\

$~~~~~~ = {\cal P} ( X ) ~\bigcup~ \{ c ~|~ c \in X \} ~\bigcup~ \{ d ~|~ d \in c \in X \} ~\bigcup~ \ldots =$ \\

$~~~~~~ = {\cal P} ( X ) ~\bigcup~ TC ( X )$ \\

11) If $x, y \in {\cal U}$ and a = \{ x, \{ y \} \}, then \\

$ TC ( a ) = \{ x, \{ y \}, y \} $ \\

since $b \in a ~\Longleftrightarrow~ b = x \,\bigvee\, b = \{ y \}$, thus $c \in b ~\Longleftrightarrow~ c = y$,
hence \\

$TC ( a ) = a ~\bigcup~ \{ c ~|~ c \in b \in a \} ~\bigcup~ \{ d ~|~ d \in c \in b \in a \} \ldots = $ \\

$~~~~~~~~ = a ~\bigcup~ \{ c ~|~ c = y \} ~\bigcup~ \{ d ~|~ d \in c = y\} \ldots = a ~\bigcup~ \{ y \} $

\hfill $\Box$ \\

The second important operation is presented in \\

{\bf Definition 3.3.} \\

We define for sets their {\it support} as follows \\

(3.3)~~~ $ SET \ni a ~~\longmapsto~~ support ( a ) = TC ( a ) \bigcap {\cal U} $ \\

Further, a set $a \in SET$ is called {\it pure}, if and only if \\

(3.4)~~~ $ support ( a ) = \phi $ \\

{\bf Note.} \\

The meaning of $support ( a )$, for a set $a \in SET$, is easy to see, based on (3.2$^*$), namely \\

(3.3$^*$)~~~ $support ( a ) = \{ x \in {\cal U} ~|~ x \in a_n \in \ldots \in a_1 \in a_0 = a,~~ n \geq 0 \}$ \\

in other words, $support ( a )$ is the set of all ur-elements $x \in {\cal U}$, if there exist any, with which {\it finite} descending sequences $x \in
a_n \in \ldots \in a_1 \in a_0 = a$, with $n \geq 0$, that start with the set $a$ do terminate.  \\

Consequently, pure sets $a \in SET$ do {\it not} have such finite descending sequences, but only {\it infinite} ones, namely \\

$ \ldots \in a_n \in a_{n - 1} \in \ldots \in a_2 \in a_1 \in a_0 = a$ \\

\hfill $\Box$ \\

Finally, the third important operation is presented in \\

{\bf Definition 3.4.} \\

(3.5)~~~ $ {\cal U} \supset A ~~\longmapsto~~ V_{afa} [ A ] = \{ a \in SET ~|~ support ( a ) \subseteq A \} $ \\

and clearly, $V_{afa} [ A ]$ is always a proper class. \\

We also denote \\

(3.6)~~~ $ V_{afa} [ \phi ] = V_{afa} = \{ a \in SET ~|~ a ~~\mbox{is a pure set} \} $

\hfill $\Box$ \\

Clearly \\

(3.7)~~~ $ V_{afa} [ A ] \subseteq SET $ \\

therefore \\

(3.8)~~~ $ V_{afa} [ A ] \bigcap {\cal U} = A \bigcap V_{afa} [ A ] = \phi,~~~ A \subseteq {\cal U} $ \\

Also, if $a \in SET$, then we have seen that \\

$TC ( a ) = \phi ~\Longleftrightarrow~ a = \phi$ \\

therefore \\

$a = \phi ~\Longrightarrow~ support ( a ) = \phi$ \\

Also, if $a \subseteq {\cal U}$, then we have seen that \\

$TC ( a ) = a$ \\

therefore \\

$support ( a ) = a$ \\

If $x, y \in {\cal U}$ and a = \{ x, \{ y \} \}, then we have seen that \\

$ TC ( a ) = \{ x, \{ y \}, y \}$ \\

thus \\

$ support ( a ) = TC ( a ) \bigcap {\cal U} = \{ x, \{ y \}, y \} \bigcap {\cal U} = \{ x, y \} $ \\ \\

{\bf 4. Graph Formulation} \\

In this section we follow the presentation in [1], without however the proofs. \\

We consider {\it directed graphs} $( N, E)$, where $N$ is the set of {\it nodes} and $V \subseteq N \times N$ is the set of {\it vertices}. A vertex
$( n, n\,' ) \in V$ can be denoted by $n \to n\,'$. Thus \\

$n_0 \to n_1 \to n_2 \to \ldots$ \\

is a finite or infinite {\it path} \\

The graph $( N, E)$ is called a {\it well-founded graph}, if and only if it has no infinite path \\

We also denote \\

$N \ni n ~\longmapsto~ [\, n > = \{\, n\,' \in N ~|~ ( n, n\,' ) \in E \,\} $ \\

$N \ni n ~\longmapsto~ < n \,] = \{\, n\,' \in N ~|~ ( n\,', n ) \in E \,\} $ \\

Given a directed graph $( N, E)$ and $n_0 \in N$, we call $( n_0 \in N, E)$ a {\it pointed graph}. \\

A pointed graph $( n_0 \in N, E)$ is called an {\it accessible graph}, if and only if \\

$~~~ \forall~ n \in N ~:~ \exists~~ path ~ n_0 \to \ldots \to n $ \\

An accessible graph $( n_0 \in N, E)$ is called a {\it tree with root} $n_0$, if and only if \\

$~~~ \forall~ n \in N ~:~ \exists\, !~~ path ~ n_0 \to \ldots \to n $ \\

A {\it decoration} of a directed graph $( N, E)$ is any mapping $S : N \longrightarrow Set$, such that \\

$~~~ \forall~ n \in N ~:~ S ( n ) = \{~ S ( n\,' ) ~|~ n\,' \in [\, n > ~\}$ \\ \\

{\bf Example 4.1.} \\

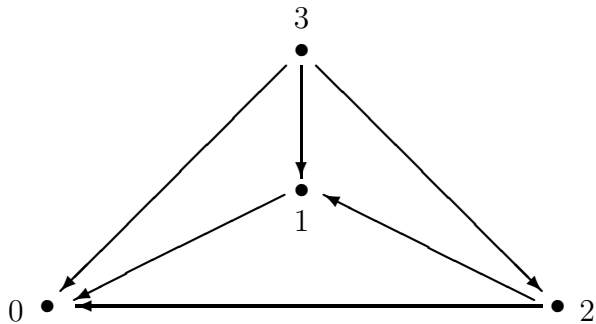
\begin{math}
\setlength{\unitlength}{0.2cm}
\thicklines
\begin{picture}(60,20)

\put(20,20){$3$}
\put(20,18){$\bullet$}
\put(19.5,17.5){\vector(-1,-1){15}}
\put(21.5,17.5){\vector(1,-1){15}}
\put(1,0.5){$0$}
\put(3.1,1){$\bullet$}
\put(39,0.5){$2$}
\put(37,1){$\bullet$}
\put(20,6.5){$1$}
\put(20,8.7){$\bullet$}
\put(36.5,1.5){\vector(-1,0){31}}
\put(19.4,8.95){\vector(-2,-1){14}}
\put(36,1.9){\vector(-2,1){14}}
\put(20.5,17.5){\vector(0,-1){7.5}}

\end{picture}
\end{math} \\

let $S$ be any decoration \\

then $[\, 0 > ~=~ \phi$, thus $S ( 0 ) = \phi = 0$ \\

and $[\, 1 > ~=~ 0$, thus $S ( 1 ) = \{ S ( 0 ) \} = \{ \phi \} = 1$ \\

while $[\, 2 > ~=~ \{ 0, 1 \}$, thus $S ( 2 ) = \{ S ( 0 ), S ( 1 ) \} = \{ 0, 1 \} = 2$ \\

finally $[\, 3 > ~=~ \{ 0, 1, 2 \}$, thus $S ( 3 ) = \{ S( 0 ) S ( 1 ), S ( 2 ) \} = \{ 0, 1, 2 \} = 3$ \\ \\

{\bf Mostowski's Collapsing Lemma 4.1.} \\

Every well-founded graph has a unique decoration.

\hfill $\Box$ \\

A {\it picture of a set} $A$ is any accessible graph $( n_0 \in N, E)$ which has a decoration $S$ such that $A = S ( n_0 )$ \\

{\bf Corollary 4.1.} \\

Every well-founded accessible graph is a picture of a unique set. \\

{\bf Proposition 4.1.} \\

Every set has a picture. \\

{\bf Proposition 4.2.} \\

The ANTI-FOUNDATION AXIOM (AFA) has the equivalent formulation : \\

$\begin{array}{|l}

\mbox{AFA AXIOM} \\ \\

\mbox{Every directed graph has a unique decoration.} \end{array} $ \\ \\

{\bf Corollary 4.2} \\

Every accessible graph is the picture of a unique set. \\

There exist non-well-founded sets. \\ \\

{\bf 5. Comments, and Beyond} \\

Self-reference has for quite a while by now happened to have acquired a rather automatic and somewhat thoughtless bad reputation as being but a source of
undesirable paradoxes. \\
One of the more memorable moments in this regard was in ancient Greece, when the man from Theba came to Athens and stated in front of Athenians that :
"All Thebans are liars !" \\
Nearer to our own days, in 1903, Russell's Paradox reformulated that ancient story within Set Theory which was then emerging as the basis of modern
mathematics, and thus further aggravated the age old negative reflexes regarding self-referentiality. \\

On the other hand, as anthropologists tell us, three fundamental themes in human thought deeply rooted in prehistoric and pre-literate times have been
self-referentiality, infinity and change. \\

Regarding the first, which is of main interest here, countless images of a snake biting its own tail are a testimony.  Also, ancient Vedic  wisdom saw it
as the foundational aspect of reality. As for the ancient Hebrews, in Exodus 3:14 of the Old Testament, they considered it to be nothing less than the
very name of God. \\
In this way, with the self-referential snake - rendered harmless as long as it is busy biting its own tail - as much as with the ancient Hindus or
Hebrews, self-referentiality was not at all a horror to be avoided by all means. On the contrary, it was a rather sacred foundational aspect of the whole
of reality ... \\

But then, later, came the man from Theba ... \\
And in our days, as a reinforcement in the very foundations of mathematics upon Set Theory, we have been facing Russell's Paradox ... \\

As it happens, however, a turn was taken in [5] back to ancient, pre-Athenian wisdom. And self-referentiality was in fact found to be of a positive
practical interest, an interest which could not be addressed in other ways, [1-3]. \\

But to return to what may be seen as more of an everyday mathematics, and actually, physics as well. Recently it was noted that such an elementary and
basic concept like {\it orthogonality} can in fact be defined in arbitrary vector spaces without any scalar product, provided a self-referential
definition is employed, [13]. \\
The relevance of that in modern physics is obvious. Indeed, in Quantum Mechanics, for instance, the standard model is based on Hilbert spaces where
orthogonality is essential and has considerable physical meaning and interpretation. \\

However, self-referentiality may turn out to have far larger and deeper impact in mathematics. In this regard, let us mention a few areas where, given
long ongoing deeper underlying difficulties so far not treatable without self-referentiality, one may at last find a more appropriate approach by
suitable self-referential definitions of basic concepts. \\

For instance, the usual concept of topology, introduced by Hausdorff in 1914, suffers among others from the fact that the respective category is not
Cartesian closed. In other words, given three topological spaces $X, Y$ and $Z$, we typically do not have the equality ${\cal C} ( X \times Y, Z ) =
{\cal C} ( X, {\cal C} ( Y, Z ) )$ between the respective spaces of continuous functions. \\
That, as well as other deficiencies of the usual concept of topological space have during the last decades been addressed by various more general
concepts of so called pseudo-topologies, [15]. The fact, however, remains that the variety of pseudo-topological concepts evolved so far, and defined of
course without self-reference, give the impression of a series of ad-hoc disparate steps which do not seem to manage to touch in a more unifying manner
the deeper meaning of topology. \\

Probability theory is another area of mathematics where the standard Kolmogorov model has manifest deficiencies. One that turns up from the very
beginning is that each point $x \in X = [ 0, 1 ]$ has the probability zero with the usual measure, thus probabilistically it is redundant. Yet the set of
all such points cannot be eliminated, since then one would remain with the empty set. And this is in sharp contradistinction with what happens in the
case of a finite or countable probability space $X$, where each point of zero probability can be eliminated, and one remains with a simplified model
$X_0$ which is isomorphic from the point of view of probability. Further well known difficulties with the standard Kolmogorv model are found in the study
of continuous time stochastic processes. \\
It is therefore an open question of some effective interest whether a self-referential definition of probability space may help in overcoming such
difficulties. In this regard it is worth noting that, while the Loeb nonstandard approach to probability brought with it a number of advantages, it has
nevertheless not been able to address satisfactorily the mentioned, as well as other difficulties. \\

The concept of computability has been of major interest during the last decades. As for the nature of its definition, the present relevance of the
related Church-Turing Thesis can be seen as showing a certain lack of sufficient insight, and thus it can appear as an inadequacy. In this regard, one
may consider the possibility of a self-referential definition of the concept of computability. \\

Complexity, among other realms in computation, is another fundamental modern concept in mathematics. And then, as its own name may possibly suggest,
perhaps, the present day simple non-self-referential definitions for it may actually be rather inappropriate ... \\

So much for avoiding the alleged horrors of self-referentiality ... \\

And as a sign of the power of persistence of age old negative connotation attached to self-referentiality, one can note that major recent contributions
to the subject still use a negative terminology, such as "non-well-founded sets" or "vicious circles" ... \\

And now, let us consider a possibly yet older, more universal, and so far incontrovertible horror, namely, that of contradiction. \\

Indeed, in this regard, there seems not to be found any controversy of any significance whatsoever throughout known human history, with all the evidence
pointing to the universal commandment of : "One must avoid contradiction !" \\

And yet, so strangely, ever since we so essentially use our modern electronic digital computers, we have been basing so much, and in such an essential
manner, on a very simple, clear and sharp contradiction. \\
Indeed, rather not consciously known to most of us, such computers - even when seen as operating only on non-negative integers - function according to
the following : \\

CONTRADICTORY SYSTEM OF AXIOMS :

\begin{itemize}

\item the well known Peano Axioms

\end{itemize}

plus the axiom :

\begin{itemize}

\item there exists M $>>$ 1, such that M + 1 = M

\end{itemize}

where the respective M, called {\it machine infinity}, may typically be larger than $10^{100}$. \\

So much for avoiding the alleged horrors of contradictions ... \\

Here however, apparently not having any known ancient wisdom to return to, a genuinely novel opening was taken with recent studies of so called {\it
inconsistent mathematics}, [11,12]. \\

What may, beyond all that, be indeed a major new opening is the development of mathematics which {\it brings together} both self-referentiality and
contradiction. And the unprecedented vastness of the respective realms that may become available in such a way is only to be guessed at present
time ... \\ \\

{\bf 6. Axioms of Set Theory} \\

For convenience, we recall here the ZFC Axioms of Set Theory, [8, p. 1]. \\

$\begin{array}{|l}

\mbox{AXIOM OF EXTENSIONALITY} \\ \\

\mbox{If two sets}~ X ~\mbox{and}~ Y ~\mbox{have the same elements, then}~ X = Y. \end{array} $ \\ \\

$\begin{array}{|l}

\mbox{AXIOM OF PAIRING} \\ \\

\mbox{For any two sets}~ a ~\mbox{and}~ b ~\mbox{there exists a set}~ \{ a, b \} ~\mbox{that contains exactly}~  a ~\mbox{and}~ b. \end{array} $ \\ \\

$\begin{array}{|l}

\mbox{AXIOM SCHEMA OF SEPARATION} \\ \\

\mbox{If}~ P ~\mbox{is a property with parameter}~ p, ~\mbox{then for any sets}~ X ~\mbox{and}~ p ~\mbox{there exists a set}~ \\ Y = \{ u \in X ~|~
P ( u, p ) \} ~\mbox{that contains all those elements}~ u \in X ~\mbox{which have property}~ P. \end{array} $ \\ \\

$\begin{array}{|l}

\mbox{AXIOM OF UNION} \\ \\

\mbox{For any set}~ X ~\mbox{there exists a set}~ Y = \bigcup X, ~\mbox{the union of all elements of}~ X. \end{array} $ \\ \\

$\begin{array}{|l}

\mbox{AXIOM OF POWER SET} \\ \\

\mbox{For any set}~ X ~\mbox{there exists a set}~ Y = {\cal P} ( X ), ~\mbox{the set of all subsets of}~ X. \end{array} $ \\ \\

$\begin{array}{|l}

\mbox{AXIOM OF INFINITY} \\ \\

\mbox{There exists an infinite set}. \end{array} $ \\ \\

$\begin{array}{|l}

\mbox{AXIOM SCHEMA OF REPLACEMENT} \\ \\

\mbox{If a class}~ F ~\mbox{is a function, then for any set}~ X ~\mbox{there exists a set}~ \\
                                           Y = F ( X ) = \{ F ( x ) ~|~ x \in X \}. \end{array} $ \\ \\

$\begin{array}{|l}

\mbox{AXIOM OF REGULARITY OR FOUNDATION} \\ \\

\mbox{Every nonempty set has an}~ \in-\,\mbox{minimal element}. \end{array} $ \\ \\

$\begin{array}{|l}

\mbox{AXIOM OF CHOICE} \\ \\

\mbox{Every family of nonempty sets has a choice function}. \end{array} $ \\ \\

Here, we also recall the additional axioms used above. \\

$\begin{array}{|l}

\mbox{ANTI-FOUNDATION AXIOM (AFA)} \\ \\

{}~~~ \forall~~ {\cal E} ~\mbox{flat system of equations} ~:~ \exists \,!~~ s~~ \mbox{solution} \end{array} $ \\ \\

An equivalent formulation of the above Axiom of Regularity or Foundation is given in \\ \\

$\begin{array}{|l}

\mbox{AXIOM OF FOUNDATION ( FA )} \\ \\

{}~~~ \forall~~ a \in SET ~:~ < a, \in \, > ~\mbox{well-founded} \end{array} $ \\ \\

$\begin{array}{|l}

\mbox{STRONG AXIOM OF PLENITUDE} \\ \\

\mbox{There is an operation}~ new ( a, b ), ~\mbox{such that} \\ \\

{}~~~ 1) ~~~\forall~~ a \in SET,~ b \subset {\cal U} ~:~ new ( a, b ) \in {\cal U} \setminus a \\ \\

{}~~~ 2) ~~~\forall~~ a, a\,' \in SET,~ a \neq a\,',~ b \subset {\cal U} ~:~ new ( a, b ) \neq new ( a\,', b ) \end{array} $ \\ \\

\end{document}